# Towards a Resolution of P = NP Conjecture

Garimella Rama Murthy

Signal Processing and Communication Research Center
International Institute of Information Technology
Hyderabad, India
rammurthy@iiit.ac.in

*Abstract*— In this research paper, the problem of optimization of a quadratic form over the convex hull generated by the corners of hypercube is attempted and solved. It is reasoned that under some conditions, the optimum occurs at the corners of hypercube. Results related to the computation of global optimum stable state ( an NP hard problem ) are discussed. An algorithm is proposed. It is hoped that the results shed light on resolving the $P \neq NP$ problem.

## I. INTRODUCTION

Optimization problems naturally arise in many fields of human endeavour. As a result, owing to the efforts of mathematicians, solution of various optimization problems is carried out using automated methods / procedures. Thus, non-experts routinely apply these automated methods. In various engineering disciplines, researchers formulated and solved constrained / unconstrained optimization problems such as linear / non-linear programming. Particularly in computer science, various discrete optimization problems ( i.e. the constraint / feasible set is a finite set ) naturally arise. For such combinatorial optimization problems, computer scientists have devised various interesting algorithms. The main criteria of efficiency of an algorithm is its computational complexity. This criteria naturally led to the class of polynomial time algorithms ( class P ) and NP complete problems ( class NP ). An interesting and difficult open research problem is whether P = NP ?

Quadratic form constitutes a second degree homogeneous polynomial in multiple variables. Optimization of a quadratic form over various constraint sets such as the bounded integer lattice arise naturally in many applications ( in science and technology ). Specifically, it is summarized in the Section 2, that optimization of quadratic form over the unit hypercube naturally arises in the case of design of Hopfield assoiciative memory. Thus, results related to such a combinatorial optimization problem find many applications.

In computer science, graphs ( directed as well as undirected ) naturally arise in many applied problems. For instance, in a directed graph, computing the minimum cut is an interesting research problem ( a special case transportation problem ). Ford-Fulkerson algorithm is the first polynomial time algorithm to compute the minimum cut in a directed graph. Several efficient polynomial time algorithms have been designed for such a problem. But, it has been realized that computing the minimum cut in an undirected graph is an NP-hard problem. So far, no polynomial time algorithm has been designed for such a problem. This research paper provides an effort in that direction.

This research paper is organized as follows. In Section 2, relevant research literature is reviewed. The author in his research efforts formulated and solved the problem of optimizing a quadratic form over the convex hull generated by the corners of unit hypercube. This result and the related ideas are documented in Section 3. In Section 4, some contributions are made towards solving the NP-hard problem of computing the global optimum stable state of a Hopfield neural network. In Section 5, the problem of optimization of quadratic / higher degree energy function on the bounded lattice is considered. The results in earlier sections are briefly generalized. In Section 6, novel distance measures on N-dimensional Eucliden space are proposed and also the structure of unit hypercube in various dimensions is explored. Finally some conclusions are reported in Section 7 .

## II. REVIEW OF RESEARCH LITERATURE

In the following, we describe how a special type of feedback neural network acts as a local optimization device with the associated quadratic energy function.

- **Contribution of Hopfield et.al:**

Hopfield neural network constitutes a discrete time nonlinear dynamical system. It is naturally associated with a weighted undirected graph G = (V,E), where V is the set of vertices and E is the set of edges. A weight value is attached to each edge and a threshold value is attached to each vertex/node of the graph. The order of the network is the number of nodes / vertices in the associated graph. Thus a discrete time Hopfield neural network of order 'N' is uniquely specified by

(A) N x N Symmetric Synaptic Weight Matrix M i.e. $M_{ij}$ denotes the weight attached to the edge from node i to node j ( node j to node i )

(B) Nx1 Threshold Vector i.e. $T_i$ denotes the threshold attached to node 'i'.

Each neuron is in one of the two states i.e. +1 or -1. Thus, the state space of such a non-linear dynamical system is the N-dimensional unit hypercube. For notational purposes, let $V_i(t)$ denote the state of node / neuron 'i' at the discrete time index 't'. Let the state of the Hopfield neural network at time 't' be denoted by the Nx1 vector $V(t)$. The state at node 'i' is updated in the following manner ( i.e. computation of next state of node 'i' )

$$V_i(t+1) = Sign\left(\sum_{j=1}^{N} M_{ij}\ V_i(t) - T_i\right) \ldots(2.1)$$

i.e. the next state at node 'i' is +1 if the term in the bracket is non-negative and -1 if it is negative. Depending on the set of nodes at which the state updation given in equation (2.1) is performed, the neural network operation is classified into the following modes:

- Serial Mode: The set of nodes at which state updation as in (2.1) is performed is exactly one i.e. at time 't' only at one of the nodes / neurons the above state updation is performed.

- Fully Parallel Mode: At time 't' , the state updation as in (2.1) is performed simultaneously at all the nodes

In the state space of discrete time Hopfiled neural network, there are certain distinguished states, called the STABLE STATES.

**Definition:** A state V(t) is called a "Stable State" if and only if
$$V(t) = Sign\ (M\ V(t) - T)\ \ldots(2.2)$$
Thus, if state dynamics of the network reaches the stable state at some time 't', it will remain there for ever i.e. no change in the state of network occurs regardless of the mode of operation of the network ( i.e. it is a fixed point in the state dynamics of discrete time Hopfield neural network ).

The following Convergence Theorem summarizes the dynamics of discrete time Hopfield neural network in the serial and parallel modes of operation. It characterizes the operation of the neural network as an associative memory.

**Theorem 1:** Let the pair N = ( M,T ) specify a Hopfield neural network. Then the following hold true:

[1] Hopfield : If N is operating in a serial mode and the elements of the diagonal of M are non-negative, the network will always converge to a stable state ( i.e. there are no cycles in the state space ).

[2] Goles: If N is operating in the fully parallel mode, the network will always converge to a stable state or to a cycle of length 2 ( i.e the cycles in the state space are of length $\leq$ 2 ).

The proof of above Theorem is based on associating the dynamics of Hopfield Neural Network (HNN) with an energy function. It is reasoned that the energy function is non-decreasing when the state of the network is updated ( at successive time instants ). Since the energy function is bounded from above, the energy will converge to some value. The next step in the proof is to show that constant energy implies a stable state is reached in the first case and atmost a cycle of length 2 is reached in the second case.

The so called energy function utilized to prove the above convergence Theorem is the following one:

$$E(t) = V^T(t)\ M\ V(t) - 2\ V^T(t)\ T\ \ldots(2.3)$$

Thus, HNN, when operating in the serial mode will always get to a stable state that corresponds to a local maximum of the energy function. Hence the Theorem suggests that Hopfield Associative Memory (HAM) could be utilized as a device for performing local/global search to compute the maximum value of the energy function.

- **Contribution of Bruck et.al**

The above Theorem implies that all the optimization problems which involve optimization of a quadratic form over the unit hypercube ( constraint / feasible set ) can be mapped to a HNN which performs a search for its optimum. One such problem is the computation of a minimum cut in an undirected graph.

For the sake of completeness, we now provide the definition of a "cut" in a graph

**Definition:** Consider a weighted, undirected graph, G with the set of edges E and the set of vertices V i.e. G = ( V, E ). Consider a subset U of V and let

$\hat{U} = V - U$. The set of edges each of which is incident at a node in U and at a node in $\hat{U}$ is called a cut in G.

**Definition:** The sum of edge weights of a cut is the weight of the cut. From among all possible cuts, the cut with minimum weight is called a Minimum cut of the graph.

In the following Theorem, proved in [BrB], the equivalence between the minimum cut and the computation of global optimum of energy function of HNN is summarized.

**Theorem 2:** Consider a Hopfield Neural Network (HNN) N= (M,T) with the thresholds at all nodes being zero i.e. $T \equiv 0$. The problem of finding the global optimum stable state ( for which the energy is maximum ) is equivalent to finding a minimum cut in the graph corresponding to N.

The author, after mastering the results in [Hop], [BrB] contemplated on removing the conditions required in Theorems 1 and 2. The fruits of such effort are documented in the following Section.

### III. OPTIMIZATION OF QUDRATIC FORMS OVER HYPERCUBE

The energy function associated with a HNN, considered in (2.3) is not exactly a quadratic form. The author questioned whether the threshold vector associated with a HNN can always be assumed to be zero ( for instance by introducing a dummy node and suitably choosing the synaptic weights from it to the other nodes ). The result of such effort is the following Lemma.

**Lemma 1**: There is no loss of generality in assuming that the threshold vector associated with a Hopfield Neural Network (HNN) is an 'all-zero' vector

Proof: Refer the argument provided in [RaN].

- Thus, it is clear that a properly chosen HNN acts as a local / global optimization device of an arbitrary quadratic form as the objective function on the constraint set being unit hypercube. Thus in the following discussion, we consider only a pure quadratic form as the energy function.

- Also, in part 1 of Theorem 1, we require the condition that the diagonal elements of symmetric synaptic weight matrix are all non-negative. We now show in the following Theorem that such a condition can be removed.

In this section, we consider the problem of maximization of quadratic form ( associated with a symmetric matrix ) over the corners of binary, symmetric hypercube. Mathematically, this set is specified precisely as follows:

$$S = \{\overline{X} = (x_1, x_2, x_3, ......, x_N) : x_i = \pm 1 \, for \, 1 \leq i \leq N\}$$
……… (3.1)

From now onwards, we call the above set simply as hypercube. This optimization problem arises in a rich class of applications. This problem is the analogue of the maximization over the hypersphere of quadratic form associated with a symmetric matrix. Rayleigh provided the solution to the optimization problem on the unit hypersphere.

A necessary condition on the optimum vector lying on the unit hypercube is now provided. This Theorem is the analogue of the maximization over the hypersphere of a quadratic form associated with a symmetric matrix. The following Theorem and other associated results were first documented in [Rama1].

**Theorem 3:** Let **B** be an arbitrary N x N real matrix. From the standpoint of maximization of the quadratic form i.e. $u^T B u$ on the hypercube, it is no loss of generality to assume that **B** is a symmetric matrix with zero diagonal elements. If $u$ maximizes the quadratic form $u^T B u$, subject to the constraint that $|u_i| \leq 1 \, for \, 1 \leq i \leq N$ ( i.e. $\overline{u}$ lies on the solid hypercube ), then

$$u = sign(Cu) \quad \text{……. (3.2)}$$

Where $C = \frac{1}{2}(B + B^T)$ with all the diagonal elements set to zero. In the above equation (i.e. eqn 3.2 ), Sign ( 0 ) is interpreted as +1 or -1 based on the requirement.

**Proof:** Any arbitrary matrix $B$ can be split into symmetric and skew-symmetric components i.e.

$$C = \frac{1}{2}(B + B^T) \, and \, \frac{1}{2}(B - B^T) \, …. (3.3)$$

Since the quadratic form associated with the skew symmetric part (matrix) is zero, as far as the optimization of quadratic form is concerned, there is no loss of generality in restricting consideration to symmetric matrices.

- It is now shown that as far as the current optimization problem is concerned, we can only consider symmetric matrices with zero diagonal elements.

Consider the quadratic form $u^T B u$, where the vector $u$ lies on the boundary of the hypercube. Since $u$ lies on the boundary, the quadratic form can be rewritten in the following form:

$$u^T C u = Trace(C) + \sum_{i=1}^{N}\sum_{j=1}^{N} u_i C_{ij} u_j$$

for $i \neq j$ …………….……..(3.4)

Since the Trace (C) is a constant, as far as the optimization over the hypercube is concerned, there is no loss of generality in restricting consideration to a matrix $\check{C}$ whose diagonal elements are all set to zero.

- In the above discussion, we assumed that the optimum of quadratic form over the convex hull of hypercube occurs on the boundary. It will be reasoned in the following discussion.

Now, we apply the discrete maximum principle [ SaW, pp.132 ] to solve the static optimization problem.

Consider a discrete time system

$Z(k+1) = u(k)$ for k=0,1, where $u(0) = u$. …… (3.5)

The criterion function to be minimized is given by

$$J^{(0)} = -\frac{1}{2} Z^T(1) \check{C} Z(1) = \theta(Z(1),1) \quad ……..(3.6)$$

The Hamiltonian is given by

$$H[Z_k, u_k, \rho_{k+1}, k] = \rho_{k+1}^T u(k) \quad ……..(3.7)$$

From the Discrete maximum principle [ SaW, pp.132 ], since $|u(0)| \leq 1$, the Hamiltonian is minimized when

$$u(0) = - \text{sign}(\rho_1) ….. (3.8)$$

From the following canonical equation [SaW, pp.133],

$$\rho_1 = \frac{\delta\theta}{\delta Z(1)} = -\check{C} Z(1) …… (3.9).$$

Thus, from (3.5), (3.8) and (3.9), we have that

$$u(0) = u = sign(\check{C} Z(1)) = sign(\check{C} u(0)) = sign(\check{C} u)$$
…………………………………………….……… (3.10).

Thus, the optimal vector $u$ satisfies the necessary condition (3.2) and it lies on the boundary of the hypercube.

Q.E.D

**Corollary:** Let **E** be an arbitrary N x N real matrix. If $u$ minimizes the quadratic form $u^T E u$, subject to the constraint

$$|u_i| = 1 \text{ for } 1 \leq i \leq n \text{ , then}$$
$$u = -sign(Cu) \quad ……………..(3.11)$$

where **C** is the symmetric matrix with zero diagonal elements obtained from **E.**

**Proof:** It may be noted that the same proof as in the above Theorem with the objective function changed from maximization to minimization of quadratic form may be used Q.E.D.

**Note:** In view of Lemma 1 and Theorem 1, a Hopfield neural network N, uniquely specified by ( W, T ) can always be assumed to be of the form where T is a zero vector and W is a symmetric synaptic weight matrix with zero diagonal elements.

Consider a weighted graph with the associated weight matrix being an arbitrary real valued matrix ( not necessarily symmetric ) S. Let the threshold value at each of the nodes be zero i.e. T≡ 0.

It is evident that
$$X^T S X = \frac{X^T (S+S^T) X}{2}.$$

Thus, from the point of view of optimization of energy function ( by the HNN ), the symmetric matrix $\frac{(S+S^T)}{2}$ can naturally be associated with the dynamics of such an associative memory.

**Definition**: The local minimum vectors of a quadratic form on the hypercube are called anti-stable states.

**Note:** It is immediate to see that if $\bar{u}$ is a stable state ( or anti-stable state ), then -$\bar{u}$ ( minus 'u' ) is also a stable state ( or anti-stable state ).

**Remark 1**: As in the case of linear programming, quadratic optimization (considered in this paper ) could be carried out using the **interior point methods** guided by the fact that global optimum over the unit hypersphere occurs at the largest eigenvector of a symmetric matrix W. The author is currently investigating this direction [Rama2], [Rama3].

**Remark 2:** The above theorem shows that optimization of a quadratic form over the convex hull generated by the corners of hypercube is equivalent to optimization just over the corners of hypercube ( i.e. Local/global optima occur only at the corners of hypercube).

**Remark 3**: The proof of the above Theorem could be given using other mathematical tools such as non-linear programming ( quadratic optimization ). Also, discrete dynamic programming based proof can be given.

**Remark 4:** It should be noted that the maximization of a quadratic form over a unit hypercube is equivalent to maximization over any hypercube. Countable union of all hypercubes is a subset of the lattice. Thus the optimum over unit hypercube could provide a decent approximation to optimization over the symmetric lattice.

**Remark 5**: Now suppose that the second sum in (3.4) does not vanish. Then, utilizing the fact that $u_i u_j = u_j u_i$, it can be rewritten as

$$\sum_{i=2}^{N}\sum_{j=1}^{i} u_i (b_{ij}+b_{ji}) u_j = u^T \widehat{B} u \quad \text{for } i \neq j \quad \ldots\ldots(3.12)$$

where $\widehat{B}$ is a lower triangular (could be upper triangular with appropriate summation) matrix with zero diagonal elements (Volterra matrix). Thus from the standpoint of the optimization over unit hypercube, it is sufficient to consider **B** to be a lower ( upper) triangular matrix with zero diagonal elements (Volterra matrix). Utilization of such a matrix could be very useful in deriving important inferences.

**Remark 6:** An upper bound on the unconstrained objective function is now given through the finite dimensional version of the Cauchy-Schwarz inequality. Let **B u = v.**

$$u^T Bu \leq \|u\| \|v\| \quad \ldots\ldots(3.13)$$

where $\|.\|$ denotes the Euclidean norm of the vector. Also equality holds if and only if

$$u = \theta v = \theta B u \ldots\ldots\ldots(3.14)$$

Thus, the result is in agreement with the Rayleigh's Theorem on optimization of quadratic form on the unit hypersphere.

- A quick argument to show that the maxima always lies on the boundary in the case of positive definite matrices is as follows:

Suppose not i.e. the extrema ( maxima ) lies inside the n-dimensional hypercube, say at $\tilde{u}$. The value of the quadratic form is given by $\tilde{u}^T B \tilde{u}$. The Euclidean norm of $\tilde{u}$ is clearly less than one. The vector $\frac{\tilde{u}}{\|\tilde{u}\|}$ which lies on the unit hypersphere gives a larger value for the quadratic form. Thus the claim is true.

**Remark 7**:

In view of equation (3.4), there is no loss of generality in assuming that the TRACE of matrix is zero for determining the stable / anti-stable states ( i.e for optimizing the quadratic form ). Since TRACE of a matrix is the sum of eigenvalues, the sum of positive valued eigenvalues is equal to the sum of negative valued eigenvalues. Hence, it is easy to see that a symmetric matrix with zero diagonal elements cannot be purely positive definite or purely negative definite. It can be assumed to be indefinite with the largest eigenvalue being a positive real number. Thus the location of stable states ( vectors ) is invariant under variation of Trace ( M ).

**Remark 8:** The stochastic versions of the problems (along the lines of Boltzmann machines) are also currently being investigated by the Author [Rama2] .

## IV. GLOBAL OPTIMUM STABLE STATE COMPUTATION

As discussed in Section II, Bruck et.al [BrB] showed that the problem of computing the maximum stable state is equivalent to that of computation of minimum cut in the associated undirected graph. This is well known to be an NP hard problem.

**Goal**: To see if a polynomial time algorithm can be discovered for the NP---hard problem of computing the minimum cut in an undirected graph. Thus we are interested in knowing whether P = NP.

In the following discussion, we consider the quadratic form associated with the matrix M ( which also can be treated as the synaptic weight matrix ).

- **Structured Class of NP-Hard Problems: Polynomial Time Algorithm:**

**Lemma 2:** Suppose the synaptic weight matrix ( connection matrix ) of the Hopfield neural network is a non-negative matrix ( i.e. all the components are non-negative real numbers ). Then, the global optimum stable state is the all ones vector i.e. [ 1 1 …. 1 ].

**Proof:** Since the energy function is a quadratic form and the variables are only allowed to assume { +1, -1 } values, the global optimum is achieved by taking the sum of all the components of the symmetric synaptic weight matrix. Hence, the vectors of all ones i.e. [ 1 1 1……1 ] is the global optimum stable state. Q.E.D.

**Lemma 3:** If a corner of unit hypercube is an eigenvector of M corresponding to positive / negative eigenvalue, then it is also a stable / anti-stable state

**Proof:** Let h be a right eigenvector of M corresponding to positive eigenvalue $\rho$. Then, we have
$$M h = \rho h \text{ with } \rho > 0.$$
Sign ( M h ) = Sign ( $\rho$ h ) = Sign ( h ) = h.
Thus h is also a stable state of M. Similar reasoning holds true when $\rho$ is a negative value.
Q.E.D.

The following remark follows from this Lemma.

**Remark 9**: Suppose we consider a vector on the hypercube ( one corner of hypercube ), say $\tilde{X}$ which is also the eigenvector of matrix M corresponding to the largest eigenvalue i.e. We have
$$M \tilde{X} = \mu_{max} \tilde{X}.$$
Then, since $\mu_{max}$ is positive, we have that

$$Sign\ (M\tilde{X}) = Sign(\mu_{max}\tilde{X}) = \tilde{X}.$$

Thus, in view of Rayleigh's Theorem, such a corner of the hypercube is also the global optimum stable state. This inference follows from the fact that the points on hypercube can be projected onto the unit hypersphere. Such a case ( where a corner of the hypercube is also the eigenvector corresponding to the maximum eigenvalue ) is very special.

- It is well known that the computation of maximum eigenvector of a symmetric matrix can be carried out using a polynomial time algorithm ( Elsner's algorithm ). Thus in such a case P = NP.

- We now provide a simple example:

**Example 1:** Consider a Hopfield neural network with two neurons. In view of Theorem 3, there is no loss of generality in assuming the synaptic weight matrix to be
$$W = \begin{bmatrix} 0 & a \\ a & 0 \end{bmatrix}, \text{ where 'a' is a real number.}$$
The eigenvalues of W are +a, -a.

Case (i): Suppose a > 0. Thus in this case the largest eigenvalue is +a.
$$W - a I = \begin{bmatrix} -a & a \\ a & -a \end{bmatrix}.$$
Hence the eigenvector corresponding to the largest eigenvalue of W can be chosen to be
$$\begin{bmatrix} 1 \\ 1 \end{bmatrix} \text{ or } \begin{bmatrix} -1 \\ -1 \end{bmatrix}.$$
It is also easy to see that, since a > 0, the global optimum stable state is
$$\begin{bmatrix} 1 \\ 1 \end{bmatrix} \text{ or } \begin{bmatrix} -1 \\ -1 \end{bmatrix}.$$
Thus, in this case, the global optimum stable state also happened to be the eigenvector corresponding to the largest eigenvalue.

Through a similar reasoning, it is easily concluded that even when a < 0, the global optimum stable state happens to be the eigenvector corresponding to the largest eigenvalue. Details are avoided for brevity.

Now, we consider the arbitrary case where the eigenvector corresponding to the largest eigenvalue is NOT a stable state

**Lemma 4:** If ' y ' is an arbitrary vector on hypercube that is projected onto the unit hypersphere and $x_0$ is the eigenvector of symmetric matrix **M** corresponding to the maximum eigen value ( on the unit hypersphere ), then we have that

$$y^T M y = \mu_{max} + 2\mu_{max}(y-x_0)^T x_0 + (y-x_0)^T M(y-x_0)$$

Proof: Let y be a vector on the hypercube that is projected onto the hypersphere. Also, let $x_0$ be the eigenvector of the symmetric synaptic weight matrix associated with the maximum eigen value. Hence the quadratic form associated with y can be expressed in the following manner.
$$y^T M y = (y - x_0 + x_0)^T M(y - x_0 + x_0) =$$
$$(y-x_0)^T M(y-x_0) + x_0^T M x_0 + 2(y-x_0)^T M x_0$$
Utilizing the fact that $x_0$ is the eigenvector corresponding to the maximum eigen value ( maximal eigenvector ) i.e. $M x_0 = \mu_{max} x_0$ and that $x_0$ lies on the unit hypersphere, that

$$y^T M y = \mu_{max} + 2\mu_{max}(y-x_0)^T x_0 + (y-x_0)^T M(y-x_0)$$

Q.E.D.

**Remark 10**: Since, by Rayleigh's theorem, it is well known that the global optimum value of a quadratic form on the unit hypersphere is the maximum eigen value i.e. $\mu_{max}$, it is clear that for all corners of the hypercube projected onto the unit hypersphere, we must necessarily have that

$$2\mu_{max}(y-x_0)^T x_0 + (y-x_0)^T M(y-x_0) \leq 0.$$

The goal is to choose a y, such that the above quantity is as less negative as possible ( so that the value of quadratic form is as close to $\mu_{max}$ as possible.

Unlike in Remark 9, suppose that
$$L = Sign(x_0) \neq Sign(Mx_0).$$
Then a natural question is to see if L can be utilized some how for arriving at the global optimum stable state. Such a question was the starting point for the following algorithm to compute the global optimum stable state.

- **Algorithm for Computation of Global Optimum Stable State of a Hopfield Neural Network:**

Step 1: Suppose the right eigenvector corresponding to the largest eigen value of M is real ( i.e. real valued components ). Compute such an eigenvector, $x_0$.

Step 2: Compute the corner, L of hypercube from $x_0$ in the following manner:
$$L = Sign(x_0).$$

Step 3: Using L as in the initial condition ( vector ), run the Hopfield neural network in the serial mode of operation.

In view of the following Lemma, eigenvector corresponding to the largest eigenvalue can always be assumed to contain real valued components.

**Lemma 5:** If A is symmetric and real matrix, then every eigenvector can be CHOSEN to contain real valued components:

**Proof:** Follows from standard arguments associated with symmetric matrices. Specifically, it is well known that any real symmetric matrix can be diagonalized by a real orthogonal matrix. Thus, one can find a system of real eigenvectors.     Q.E.D.

- In view of Lemma 2, consider the class of synaptic weight matrices ( connection matrices ) of Hopfield neural networks that are "irreducible" non-negative matrices ( i.e. necessarily all the components are non-negative real numbers ).

For such non-negative matrices, the following Theorem enables us to understand the structure of eigenvector corresponding to the largest eigenvalue.

**Theorem 4 ( Perron—Frobenius Theorem ):**
Let P be an irreducible non-negative matrix. Then P has an eigenvalue $\alpha$ which is real, positive, and simple. For any other eigenvalue $\delta$ of P we have $|\delta| \leq \alpha$. To this maximal eigenvalue $\alpha$ there corresponds a strictly positive eigenvector, say $f$.

Thus, we have
$$Pf = \alpha f \text{ with } \alpha > 0.$$
Let
$$Sign(f) = L.$$
Thus, for such class of symmetric non-negative synaptic weight matrices,
 SIGNUM ( MAXIMAL EIGENVECTOR)
happens to be the global optimum stable state.

In view of Lemma 4 and above remark 10, the claim is that the global optimum stable state is reached through the above procedure. The formal proof that global optimum stable state is reached is provided below.

**Theorem 5:** The above algorithm converges to the global optimum stable state using the vector L as the initial condition.

**Proof:** In view of the results in [1] ( [BrB] ), the idea is to reason that the vector L is in the domain of attraction of the global optimum stable state. Equivalently using the results in [1], we want to show that the initial vector is in the coding sphere of the codeword corresponding to the global optimum stable state.

Let $y_0$ be the global optimum stable state / vector on the unit hypercube. Let $K_0$ be one among the other stable states. Thus, the idea is to reason that the Hamming distance between 'L' and $y_0$ i.e. $d_H(L, y_0)$ is smaller than the Hamming distance between L and $K_0$ i.e. $d_H(L, K_0)$ i.e.

To reason that
$$d_H(L, y_0) < d_H(L, K_0).$$

The proof is by contradiction i.e. say

$$d_H(L, K_0) < d_H(L, y_0) \quad ....(4.1)$$

We know that the "sign" structure of the vectors 'L' and $x_0$ is exactly the same. More explicitly all the components of 'L' and $x_0$ have the same sign ( positive or negative ).

Since the three vectors $\{L, y_0, K_0\}$ lie on the unit hypercube, we consider various possibilities with respect to the "sign" structure of those vectors. Thus, we define the

following sets:

$\bar{A}$.....Set of components of vectors $\{y_0$ and $K_0\}$ that agree ( both of them ) in "sign" with those of the vector $x_0$ ( and hence 'L' ).

$\bar{B}$.....Set of components of vectors $\{y_0$ and $K_0\}$ that DONOT agree ( both of them ) in "sign" with those of the vector $x_0$ ( and hence 'L' ).

$\bar{C}$.....Set of components of vector $y_0$ where only $y_0$ differs in "sign" from those components of vector $x_0$.

$\bar{D}$.....Set of components of vector $K_0$ where only $K_0$ differs in "sign" from those components of vector $x_0$.

By the hypothesis, the cardinality of set $\bar{C}$, i.e. $|\bar{C}|$ is atleast one larger than the cardinality of the set $\bar{D}$ i.e. $|\bar{D}|$. For concreteness, we first consider the case where $|\bar{C}| = |\bar{D}| + 1$.

To illustrate the argument, we first consider the case where only the last component of $y_0$ differs from that of $x_0$ in sign ( but not $K_0$ ) and all other components ( of both $y_0$, $K_0$ ) either agree or disagree in sign with those of $x_0$.

To proceed with the proof argument, the vectors L, $y_0$, $K_0$ ( lying on the unit hypercube ) are projected onto the unit hypersphere through the following transformation: Let the projected vectors be Q, R i.e.

$$Q = \frac{y_0}{\sqrt{N}}, \qquad R = \frac{K_0}{\sqrt{N}},$$

where N is the dimension of the symmetric matrix M..

Thus, we want to reason ( if the Hamming distance condition specified above i.e. equation (4.1) is satisfied i.e. our hypothesis ) that the value of quadratic form associated with the vectors Q, R satisfies the following inequality i.e.

$$R^T M\ R\ >\ Q^T M\ Q.$$

Hence the idea is to thus arrive at a contradiction to the fact that $y_0$ is the global optimum stable state.

In view of Lemma 4, we have the following expressions for $Q^T M Q$ and $R^T M R$ :

$$Q^T M\ Q\ =\ \mu_{max} + 2\mu_{max}(Q - x_0)^T x_0 + (Q - x_0)^T M (Q - x_0)$$

$$R^T M\ R\ =\ \mu_{max} + 2\mu_{max}(R - x_0)^T x_0 + (R - x_0)^T M (R - x_0)$$

Equivalently, we effectively want to show that
$$2\mu_{max}(R - Q)^T x_0 + (R - x_0)^T M (R - x_0)$$
$$-\ (Q - x_0)^T M (Q - x_0) > 0.$$
Let us label the terms as follows:

$$(\text{I}) = 2\mu_{max}(R - Q)^T x_0$$

$$(\text{II}) = (R - x_0)^T M (R - x_0) - (Q - x_0)^T M (Q - x_0) > 0.$$

We want to show that
$$(\text{I}) + (\text{II}) > 0.$$

To prove such an inequality, we first partition the vectors Q, $x_0$, R ( lying on the unit hypersphere ) into two parts:

- (i) part (A) where components of Q, R simultaneously agree or disagree in sign with those of $x_0$,

- (ii) part (B) : Last component of Q, that disagrees in sign with the last component of $x_0$. But the last component of R agrees in sign with that of $x_0$.

Thus, the vectors $x_0$, Q, R are given by
$$x_0 = \begin{bmatrix} x_0^A \\ x_0^B \end{bmatrix}, \ Q = \begin{bmatrix} Q_A \\ Q_B \end{bmatrix}, \ R = \begin{bmatrix} R_A \\ R_B \end{bmatrix}, \text{ where}$$
$Q_B$, $R_B$ are scalars. Also the components of $\{Q_A, R_A\}$ are simulataneously either $+\frac{1}{\sqrt{N}}$ or $-\frac{1}{\sqrt{N}}$. Thus except for the last component, all other components of the vector R-Q are all zeroes. Further suppose that the last component of $x_0$ i.e. $x_0^B$ is $-\theta$, with $\theta > 0$. Then it is easy to see that
$$(R_B - Q_B) = -\frac{2}{\sqrt{N}}.$$

Summarizing
$$(R - Q)^T = \begin{bmatrix} 0\ 0\ ...\ 0\ \frac{-2}{\sqrt{N}} \end{bmatrix}.$$
Hence, we have that
$$(R - Q)^T x_0\ =\ \frac{2}{\sqrt{N}}\ \theta.$$

Thus $(\text{I}) = \frac{4\ \mu_{max}\ \theta}{\sqrt{N}}$

which is strictly greater than zero. Similarly, even when the last component of $x_0$ is $+\theta$ with $\theta > 0$, it is easy to reason that $(R - Q)^T x_0$ is strictly greater than zero.

Now we consider the other term i.e Term ( II ):
We first partition M into a block structured matrix i.e.

$$M = \begin{bmatrix} M_{11} & M_{12} \\ M_{21} & \beta \end{bmatrix}, \text{ where } \beta \text{ is a scalar and}$$
$$M_{21}\ =\ M_{12}^T.$$
We also partition the vectors $(Q - x_0)$, $(R - x_0)$ in the following manner:

$$(Q - x_0) = \begin{bmatrix} F^{(1)} \\ G^{(1)} \end{bmatrix} \; ; \; (R - x_0) = \begin{bmatrix} F^{(2)} \\ G^{(2)} \end{bmatrix},$$

Where $G^{(1)}, G^{(2)}$ are scalars. As per the partitioning procedure, it is clear that
$$F^{(1)} = F^{(2)}.$$
Also, let us consider the case where the last component of $x_0$ i.e. $x_0^{(B)}$ is $-\theta$, with $\theta > 0$. In such a case
$$G^{(1)} = \frac{1}{\sqrt{N}} + \theta \; ; \; G^{(2)} = -\frac{1}{\sqrt{N}} + \theta.$$

**Note:** For the case, where the last component of $x_0$ is $+\theta$ with $\theta > 0$, all the following equations are suitably modified. Details are avoided for brevity.

In term (II), the following definition are utilized:
$$H = (R - x_0)^T M (R - x_0) \text{ and }$$
$$J = (Q - x_0)^T M (Q - x_0).$$

Thus  (II) = H – J.

In view of partitioning of the matrix M and vectors $(Q - x_0)$, $(R - x_0)$; we have that

$$H = F^{(2)^T} M_{11} F^{(2)} + 2 F^{(2)^T} M_{12} G^{(2)} + G^{(2)^T} M_{22} G^{(2)}$$

$$J = F^{(1)^T} M_{11} F^{(1)} + 2 F^{(1)^T} M_{12} G^{(1)} + G^{(1)^T} M_{22} G^{(1)}.$$

Using the fact that $F^{(1)} = F^{(2)}$, we have that

$$H - J = 2 F^{(1)^T} M_{12} \left( G^{(2)} - G^{(1)} \right) + \beta \left[ (G^{(2)})^2 - (G^{(1)})^2 \right]$$

Let $F^{(1)^T} M_{12} = M_{12}^T F^{(1)} = \gamma$. Thus, we have that

$$H - J = 2\gamma \left( \theta - \frac{1}{\sqrt{N}} - \theta - \frac{1}{\sqrt{N}} \right) +$$
$$\beta \left[ (\theta - \frac{1}{\sqrt{N}})^2 - (\theta + \frac{1}{\sqrt{N}})^2 \right]$$

$$= -\frac{4\gamma}{\sqrt{N}} + \beta \left( \frac{-4\theta}{\sqrt{N}} \right).$$

Hence, we have the following expression for (I)+(II)

$$(I)+(II) = \frac{4\mu_{max}\theta}{\sqrt{N}} - -\frac{4\gamma}{\sqrt{N}} + \frac{-4\theta\beta}{\sqrt{N}}.$$

$$= \frac{4}{\sqrt{N}} \left[ -\gamma - (\beta\theta - \mu_{max}\theta) \right].$$

But since $x_0$ is the eigenvector of M corresponding to the largest eigenvalue $\mu_{max}$, we have that

$$M_{12}^T x_0^{(A)} - \beta\theta = -\mu_{max}\theta.$$

$$M_{12}^T x_0^{(A)} = \beta\theta - \mu_{max}\theta.$$

Hence we necessarily have that

$$\begin{aligned}
(I) + (II) &= \frac{4}{\sqrt{N}} \left[ -\gamma - M_{12}^T x_0^{(A)} \right] \\
&= \frac{4}{\sqrt{N}} \left[ -M_{12}^T F^{(1)} - M_{12}^T x_0^{(A)} \right] \\
&= \frac{4}{\sqrt{N}} \left[ -M_{12}^T (Q^{(A)} - x_0^{(A)}) - M_{12}^T x_0^{(A)} \right] \\
&= \frac{4}{\sqrt{N}} \left[ -M_{12}^T Q^{(A)} \right].
\end{aligned}$$

We first note that $M_{12}$ is constrained by the fact that $y_0$ is a stable state. Thus

$$Sign(M y_0) = y_0.$$

Or equivalently

$$Sign\left(M \frac{y_0}{\sqrt{N}}\right) = y_0 = Sign(M Q).$$

Thus, we necessarily have

$$Sign\left( \begin{bmatrix} M_{11} & M_{12} \\ M_{12}^T & \beta \end{bmatrix} \begin{bmatrix} Q^{(A)} \\ \frac{1}{\sqrt{N}} \end{bmatrix} \right) = \begin{pmatrix} y_0^{(A)} \\ +1 \end{pmatrix}.$$

From Theorem (3), we have freedom in choosing the diagonal elements of M (since Trace(M) contributes a constant value to the value of quadratic form on the hypercube). Thus by a suitable choice of $\beta$, we can ensure that
$$M_{12}^T Q^{(A)} < 0.$$
It can be easily reasoned that the freedom in choosing the diagonal elements of M can be capitalized to ensure that

$$Sign\left( M_{11} Q^{(A)} + \frac{1}{\sqrt{N}} M_{12} \right) = Y_0^{(A)}$$

is always satisfied ( by a proper choice of diagonal elements of M ).

Thus, we arrive at the desired contradiction to the fact that $y_0$ is a global optimum stable state. ( and $K_0$ is not the one ). Thus, the vector L is in the domain of attraction of the global optimum stable state. Hence, with this choice of initial condition, when the Hopfield neural network is run in the serial mode, global optimum stable state is reached.

Similar argument ( with block matrices ) can easily be provided for the case where $|\bar{C}| \geq |\bar{D}| + 2$. Detailed argument is avoided for brevity.

Q.E.D.

- **Calculation of Computational Complexity associated with the Above Algorithm:**

The above algorithm involves the following computations:

(A) Computation of eigenvector corresponding to the largest eigenvalue of the symmetric matrix ( connection matrix of Hopfield neural network ). A polynomial time algorithm for such a task is already available.

(B) Using the associated vector as the initial condition, running the Hopfield neural network in serial mode until global optimum stable state is reached. It is possible to bound the number of computations for this task.

**Interesting Generalizations:**
In [RaN], the authors proposed a generalization of Hopfield neural network, called the Complex Amari-Hopfield neural network. In this case the "complex hypercube" constitutes the state space of the neural network. The synaptic weight matrix of the network is a Hermitian symmetric matrix.

Thus, in the case of such complex Amari-Hopfield network, the "complex signum" function of the largest eigenvector ( i.e. eigenvector corresponding to the largest eigenvalue ) is utilized as the initial condition to run the network in "serial mode". It is reasoned that the global optimum stable state is reached through such a procedure.

- **Hopfield Neural Network: Associated One Step Associative Memory: Efficient Algorithm for Minimum Cut Computation:**
We now propose a method which reduces the computational complexity of the method of computing the global optimum stable state [Rama3]. The author is actively pursuing the problem of minimum cut computation in an undirected graph [Rama2].

**Lemma 6:** Given a linear block code, a neural network can be constructed in such a way that every local maximum of the energy function corresponds to a codeword and every codeword corresponds to a local maximum.

**Proof:** Refer the paper by Bruck et.al [BrB].

It has been shown in [BrB] that a graph theoretic code is naturally associated with a Hopfield network ( with the associated quadratic energy function ). The local and global optima of the energy function are the code words**.**

**Goal: To compute the global optimum stable state ( i.e. global optimum of the energy function ) using the associated graph theoretic encoder.**

To achieve the goal, once again the largest real eigenvector is utilized as the basis for determining the information word that will be mapped to a global optimum stable state/ codeword ( using the associated graph theoretic encoder ).

- **Tighter Lower Bound on the Spectral Radius**:
We discussed the idea of projecting points, X on the unit hypercube onto the unit hypersphere, Y using the following transformation:
$$Y = \frac{X}{||X||} = \frac{X}{\sqrt{N}},$$
where $L^2 - norm$ is utilized.
Hence, we have that
$$Y^T M Y = \frac{1}{N} X^T M X.$$
Using the Rayleigh's Theorem, we have the following inequality satisfied by the eigenvalues of M
$$\mu_{min} \leq Y^T M Y \leq \mu_{max}.$$
Hence, we have that
$$N \mu_{min} \leq X^T M X \leq N \mu_{max}.$$

We have reasoned in Theorem 3 that
$$X^T M X = \frac{1}{2} X^T (M + M^T) X = X^T W X$$

It is obvious that if M is a non-negative matrix, then W is also a non-negative matrix. More interestingly, even if "M" is not a non-negative matrix, W can be a non-negative matrix. This can happen if
$$M_{ij} + M_{ji} \geq 0 \quad for\ all\ i,j.$$

In the following discussion, we consider matrices M for which the corresponding matrix W is a non-negative matrix.

Since "X" lies on the unit hypercube and W is non-negatifve, the all-ones vector i.e. $[1\ 1\ ....1]^T$ maximizes the quadratic form ( as discussed in Lemma 2 ). Thus the stable value corresponding to such stable vector is
$$\sum_{i=1}^{N} \sum_{j=1}^{N} W_{ij}.$$
Thus, we have that
$$\mu_{max} \geq \frac{1}{N} \left[ \sum_{i=1}^{N} \sum_{j=1}^{N} W_{ij} \right].$$

Now, we argue that this bound is tighter than the related well known lower bound derived earlier. From Linear algebra, we know that
$$\tau \leq \mu_{max} \leq \delta,  \quad \text{where}$$
$\tau, \delta$ are the minimum and maximum row sums respectively ( more precisely the sums of absolute value of the elements of rows ).
It is immediate that
$$\frac{1}{N}\left[\sum_{i=1}^{N}\sum_{j=1}^{N} W_{ij}\right] \geq \tau$$
Hence, we have a tighter lower bound on the spectral radius of W ( associated with the matrix M ).

**Remark 11**: Using a similar argument, it can easily be shown that the same lower bound can be derived when the constraint set is the symmetric bounded lattice. Detailed derivation for this more general case can be found in [Rama2]. Complete derivation is avoided for brevity.

**Remark 12:**
It is easy to verify that
$$\frac{1}{N}\left[\sum_{i=1}^{N}\sum_{j=1}^{N} W_{ij}\right]$$
is a matrix norm on the class of N x N non-negative matrices. The properties and results associated with such a norm are currently being investigated [Rama2].

- **Lebesgue Decomposition : Optimization of Quadratic Forms:**
As discussed in [Rama5], any finite dimensional linear opertator, represented by a matrix, B can be decomposed as
$$B = B^{(+)} - B^{(-)},$$
where $B^{(+)}$ ( positive part of B ) contains the non-negative elements of B in the same position and $B^{(-)}$ ( negative part of B ) contains the non-positive elements of B in the same position. It is immediate that the quadratic form associated with B can be decomposed in the following manner:
$$X^T B X = X^T B^{(+)} X - X^T B^{(-)} X.$$
Now define
$$C^{(+)} = \frac{B^{(+)} + [B^{(+)}]^T}{2} \text{ and } C^{(-)} = \frac{B^{(-)} + [B^{(-)}]^T}{2}$$
Thus, we have
$$X^T B X = X^T C^{(+)} X - X^T C^{(-)} X,$$
where $C^{(+)}$ and $C^{(-)}$ are symmetric matrices

## V OPTIMIZATION ( LOCAL/GLOBAL MAXIMIZATION / MINIMIZATION ) OF QUADRATIC / HIGHER DEGREE ENERGY FUNCTION ON THE BOUNDED LATTICE: HIGHER DIMENSIONAL GENERALIZATIONS:

- There are many discrete / combinatorial optimization problems dealing with optimization of quadratic energy function on the bounded integer lattice ( i.e. finitely many feasible points ). Optimization of a quadratic form over the unit hypercube is a special case problem. As in the above discussion, the fact that eigenvector corresponding to the maximum eigenvalue provides the global optimum value ( of quadratic form ) on the unit hypersphere is capitalized.

**Approach for Optimization of Quadratic form on the symmetric bounded lattice:**

- Suppose the feasible points on the bounded lattice are projected onto the unit hypersphere. The lattice point nearest to the maximum eigenvector is determined ( Closest Lattice Point Problem ) in the following manner:

Quantization of Components of Eigenvector Corresponding to the Largest Eigenvalue:
Suppose $q_o$ is a component of the eigenvector Q corresponding to the largest eigenvalue. Then using the lower / upper ceiling function, it is quantized. In otherwords, in the quantized vector of Q, the corresponding component takes the following values:

- Lower Ceiling ( $q_0$ )… $q_0$ is truncated to nearest integer

- Upper Ceiling ( $q_0$ )… $q_0$ is rounded off to nearest integer

Using the quantized vector, the global optimum solution is determined using the results in [1] i.e. [BrB].

In the following, we consider one possible higher dimensional generalization:

**Optimization of Quadratic Form on Higher Dimensional Unit Hypercube**:

- In [7] i.e. [Rama4], the author conceived the idea of multi-dimensional Hopfield neural

network in which the state of the non-linear dynamical system is a higher dimensional array i.e. a tensor. The connection structure is captured using a fully symmetric tensor $\widehat{M}$. Such higher dimensional Hopfield associative memory also optimizes a quadratic energy function. Associated with such a multi-dimensional neural network, an interesting convergence theorem is proved ( in the same spirit of the one-dimensional case ). In view of the results in the previous section, the eigentensor of $\widehat{M}$ associated with the largest eigenvalue, say $\widehat{Z}$, is computed. Define a tensor associated with $\widehat{Z}$ as
$$\widehat{L} = Sign\ (\ \widehat{Z}\ )$$
i.e the components of the tensor $\widehat{L}$ are sign ( +1 or -1 ) of the components of the tensor $\widehat{Z}$. It should be noted that Sign ( zero ) is consistently taken as +1 or -1. Using the tensor $\widehat{L}$ lying on the multi-dimensional hypercube as the initial condition, the multi-dimensional Hopfield network is run in the serial mode. As in the previous section, the claim is that the associative memory running in serial mode converges to the global optimum stable state tensor.

- Consider a homogeneous multi-variate polynomial of EVEN degree ( i.e. higher degree than quadratic forms ) . .It is well known that such a multi-variate polynomial can naturally be specified using a Tensor linear operator. Let us specifically consider multi-variate homogeneous form of EVEN degree ( strictly larger than 2 ) generated using a symmetric Tensor linear operator. To optimize such a higher degree energy function, over the higher dimensional hypersphere, we compute the eigentensor corresponding to the largest eigenvalue. Quantizing the eigentensor corresponding to the largest eigenvalue, we arrive at a tensor that can be mapped possibly to the global optimum solution ( on generalizing the results in [BrB] ) on the higher dimensional symmetric bounded lattice.

- Suppose, we consider the generator matrix G associated with a linear block code. Using the results in [BrB], the polynomial representation of G, also called the "energy function" is easily determined. Suppose such an energy function is a multi-variate homogeneous polynomial ( also called multi-linear form in the literature ). It can be expressed in terms of an associated tensor W. Compute the eigentensor of W corresponding to the largest eigenvalue. Quantizing such an eigentensor ( as discussed previously ) will lead to a WORD ( tensor ) in the coding sphere associated with the codeword tensor that is the global optimum of the energy function.

- The following generalization of Rayleigh's Theorem is expected.
 The problem is to optimize an even degree ( > 2 ) homogeneous multivariate form captured by a higher order tensor ( than a matrix ) ( i.e. a multi-linear form) on the multi-dimensional unit hypersphere. It is expected that the eigentensors are the local optimum tensors of such a form ( with the corresponding eigenvalues ).
 With such a generalization, the results of this research paper are generalized to higher dimensions. Details can be found in [Rama2].
 The results and concepts utilized in section IV naturally motivated us to understand the structure of unit hypercube in various dimensions. Also, using Theorem 3, Lemma 3, we generalized the idea of Hopfield for the synthesis of Associative Memory ( i.e. motivated by the Hopfield's synthesis of Associative Memory ). Also a novel distance measure is introduced in the following Section.

## VI NOVEL DISTANCE MEASURE ON N-D EUCLIDEAN SPACE: STRUCTURE OF UNIT HYPERCUBE in VARIOUS DIMENSIONS:

The idea utilized in the above algorithm ( for computing the global optimum stable state of a Hopfiled neural network ) naturally motivated us to design and study a "distance measure" between any two points on the N-dimensional Euclidean space. We need the following preliminary discussion.
 Consider { +1, -1 } vectors on the unit hypercube. Consider the mapping:
$$+1\ \rightarrow\ +1$$
$$-1\ \rightarrow\ 0$$

**Definition:** "Hamming—Like" distance between any two vectors on the unit hypercube is the Hamming distance between them ( with the above mapping )

**Remark 13:** Consider any two vectors P, Q on the N-dimensional Euclidean space. By dividing each of the components of P, Q by the corresponding

Euclidean norm, we arrive at vectors on the unit hypersphere ( i.e. after normalization by the corresponding Euclidean norm ).

Using this approach, we now restrict consideration to vectors on the unit hypersphere obtained from the corresponding vectors on the N-d Euclidean space.

**Definition:** Consider any two points X, Y on the unit hypersphere. Define

$$Z = \text{Sign}(X)$$

$$W = \text{Sign}(Y)$$

i.e. the components of vector Z are obtained as the sign of the corresponding components of X. It should be noted that the sign of ZERO component is consistently defined as +1 or -1 ( i.e. Sign(0) = +1 or -1 ).

The "induced Hamming-type Distance" between { X, Y } is defined as the Hamming like distance between Z, W.

The above discussion is now generalized to determine the " generalized induced hamming distance" between two vectors on the N-dimensional Euclidean space.

We consider two points / vectors $X_1$, $X_2$ on the "bounded" ( in the sensor of Euclidean norm ) N-dimensional Euclidean space. We quantize the components of the two vectors using the ceiling function

$Y_1(j) = \text{Ceiling}[X_1(j)]$ for $1 \leq j \leq N$ and

$Y_2(j) = \text{Ceiling}[X_2(j)]$ for $1 \leq j \leq N$.

**Note:** In the above equations, we can use Lower Ceiling or Upper Ceiling function. Hence, in effect the components of vectors are rounded off / truncated to the nearest integer. It should be noted that after quantization, the components of the vectors will be positive or negative integers. Boundedness of vectors ensures that after quantization, all the components are below certain integer. Thus, the operation of quantization ensures that the vectors lie on the bounded lattice.

**Definition:** Consider any two bounded vectors $X_1$, $X_2$ lying on the N-dimensional Eucliden Space. Let
$Y_1 = \text{Ceiling}(X_1)$
$Y_2 = \text{Ceiling}(X_2)$.
The "generalized induced Hamming Distance" between the bounded vectors $X_1$, $X_2$ is defined as the Hamming distance between the vectors $Y_1$, $Y_2$.

**Note**: We can also define "induced Manhattan distance" between $X_1$, $X_2$ as the Manhattan distance between the vectors $Y_1$, $Y_2$.

In the above two definitions, a countable / finite set of vectors / points is extracted ( through the process of appropriate quantization ) from the uncountable set of vectors. Using Hamming / Manhattan distance on the countable / finite set, the "induced" Hamming / Manhattan distance is defined.

Using the above definition, we derive the following interesting Lemma. The lemma requires the following well known definition:

**Definition:** Two vectors X, Y lying on the unit hypercube are orthogonal if their inner product is zero i.e.

$$<X, Y> = 0.$$

**Lemma 7**: If the dimension of hypercube is odd, then there are NO orthogonal vectors lying on the hypercube.

**Proof**: Consider two vectors X, Y lying on the unit hypercube. Let the "Hamming—Like distance" between them be 'd'. The inner product of X, Y is given by'

$<X, Y> = \{$ Number of components where X, Y agree $\} - \{$ Number of components where X, Y disagree $\}$

$= (N - d) - d = N - 2d$.

Thus for the two vectors X, Y to be orthogonal, it is necessary and sufficient that "N" is an even number.
**Q.E.D.**

Thus, we can study the properties of the function $f(d) = N - 2d$ and interpret it suitably. Suppose, we fix "X" and vary the vector Y ( among other vectors on the unit hypercube. ). We have that

Number of vectors at a Hamming Like distance of 'k' to the vector $X = \binom{N}{k}$ for $0 \leq k \leq N$.

Thus there are $\binom{N}{k}$ vectors for which $<X, Y> = N - 2k$.

The above lemma sheds light on the structure of unit hypercube when the dimension of it is even / odd. The generalization of above Lemma to higher dimensional { +1, -1 } valued tensors ( i.e. inner product of two { +1, -1 } valued tensors ) is straightforward and is avoded for brevity.

- **On The Existence of Associative Memory Synthesized by Hopfield:**

Hopfield synthesized a real valued synaptic weight matrix from the patterns to be stored in such a way that the network so obtained has these patterns as stable states. The weight matrix given by Hopfield is as follows:

$$W = \sum_{j=1}^{S}(X_j X_j^T - I)$$

where $S$ is the number of patterns to be stored ($S < N$), $I$ is the identity matrix and $X_1, X_2 \ldots X_S$ are the orthogonal real patterns (lying on the real hypercube) to be stored. Thus it is easy to see that
$$W X_k = (N - S) X_k \text{ with } S < N.$$
Hence, the corner of hypecube $X_k$ is also an eigenvector corresponding to the positive eigenvalue ($N - S$). Thus the spectral representation of the connection matrix of associative memory synthesized by Hopfield is given by

$$W = \sum_{j=1}^{S}(N-s)\frac{X_j}{\sqrt{N}}\frac{X_j^T}{\sqrt{N}} = \sum_{j=1}^{S}\frac{(N-s)}{N} X_j X_j^T$$

Thus, in view of the above Lemma 7, we have the following result.

**Lemma 8:** Hopfield construction of associative memory exists only when the dimension of the hypercube is EVEN.

Proof: It follows from Lemma 5

In view of the fact that the synaptic weight matrix can be an arbitrary symmetric matrix, we can generalize the Hopfield construction (of associative memory) in the following manner:

Let $\{\mu_j\}_{j=1}^{S}$ be desired positive eigenvalues (with $N\mu_j$ being the desired stable value) and let $\{X_j\}_{j=1}^{S}$ be the desired stable states. Then it is easy to see that the following symmetric matrix constitutes the desired synaptic weight matrix of Hopfield neural network:

$$W = \sum_{j=1}^{S} \mu_j X_j X_j^T$$

Once again, in this case N must be even.

In the same spirit, we now synthesize a synaptic weight matrix, with desired stable/anti-stable values and the corresponding stable / anti-stable states. Let $\{X_j\}_{j=1}^{S}$ be desired orthogonal stable states and $\{Y_j\}_{j=1}^{L}$ be the desired orthogonal anti-stable states. These "L" vectors are mutually orthogonal. Let the desired stable states be eigenvectors corresponding to positive eigenvalues and let the desired anti-stable states be eigenvectors corresponding to negative eigenvalues. The spectral representation of desired synaptice weight matrix is given by

$$W = \sum_{j=1}^{S}\frac{\mu_j}{N} X_j X_j^T - \sum_{j=1}^{L}\frac{\beta_j}{N} Y_j Y_j^T$$

where $\mu_j's$ are desired positive eigenvalues and $-\beta_j's$ are desired negative eigenvalues. In view of Ramark 7, the sum of positive and negative eigenvalues is equal to zero.

Hence the above construction provides a method of arriving at desired energy landscape (with orthogonal stable /anti-stable states and the corresponding positive / negative energy values).

In view of Lemma 7, when the dimension of the unit hypercube is odd, there are only the following two possibilities:

CASE A: Only one corner of the hypercube is an eigenvector in the spectral representation of W

CASE B: None of the corners of the hypercube is an eigenvector in the spectral representation of W.

- In case A, if the corresponding eigenvalue is the maximum eigenvalue of W, then it is also the global optimum stable state

- In case B, if the spectral representation of W is of the following form (Rank one matrix)

$$W = \gamma\ f_i f_i^T \text{ with all the other eignvalues are zeroes, then}$$

$$F = \text{Sign}(f_i)$$

is the global optimum stable state. The logical reasoning is fairly simple and avoided for brevity.

- In the spirit of similar idea in coding theory, we conceived the following definition

**Definition:** The Weight of a Corner on the hypercube is defined to be the number of 1's (ones) in the associated { +1, -1 } vector.

We are interested in arriving at the weight distribution i.e. Number of { +1, -1 } vectors on the unit hypercube with a certain weight i.e. G(k)

G(k) : Number of corners of N-dimensional unit hypercube with weight k

It is easy to see that
$$G(k) = \binom{N}{k} \text{ for } 0 \leq k \leq N.$$

On normalizing ( dividing ) $G(k)\ by\ 2^N$, we arrive at a binomial probability mass function with success rate $\frac{1}{2}$ i.e. B ( N, $\frac{1}{2}$ ) i.e.

$$\hat{G}(k) = \frac{1}{2^N}\binom{N}{k} = \binom{N}{k}\frac{1}{2^k}\frac{1}{2^{N-k}}$$

Note: Deeper implications of Sylvester's law of inertia are explored {Rama2].

## VII CONCLUSION

In this research paper, it is shown that optimizing the quadratic form over the convex hull generated by the corners of hypercube is equivalent to optimization over just the corners of hypercube. Some results related to the computation of global optimum stable state are discussed.

## ACKNOWLEDGMENT

The author thanks Alexander Seliverstov of Kharkevich Institute, Russia for comments on Theorem 1.

## REFERENCES

[1] [BrB] J. Bruck and M. Blaum,"Neural Networks, Error Correcting Codes and Polynomials over the Binary Cube," IEEE Transactions on Information Theory, Vol.35, No.5, September 1989.

[2] [Hop] J.J. Hopfield, "Neural Networks and Physical Systems with Emergent Collective Computational Abilities," Proceedings of National Academy Of Sciences, USA Vol. 79, pp. 2554-2558, 1982

[3] [Rama1]G. Rama Murthy,"Optimal Signal Design for Magnetic and Optical Recording Channels, " Bellcore Technical Memorandum, TM-NWT-018026, April 1st , 1991

[4] [Rama2]G. Rama Murthy,"Efficient Algorithms for Computation of Minimum Cut in an Undirected Graph,"Manuscript in Preparation

[5] [Rama3]G. Rama Murthy,"Multi-dimensional Neural Networks : Unified Theory, " Research monograph published by NEW AGE INTERNATIONAL PUBLISHERS, NEW DELHI, 2007.

[6] [RaN] G. Rama Murthy and B. Nischal," Hopfield-Amari Neural Network : Minimization of Quadratic forms," The 6th International Conference on Soft Computing and Intelligent Systems, Kobe Convention Center (Kobe Portopia Hotel) November 20-24, 2012, Kobe, Japan.

[7] [Rama4] G. Rama Murthy,"Multi/Infinite Dimensional Neural Networks: Multi / Infinite Dimensional Logic Theory," International Journal of Neural Networks, Volume 15, No. 3, June 2005

[8] [SaW]A.P. Sage and C.C. White,"Optimum Systems Control," Prentice Hall Inc, 1977

[9] [Rama5] G. Rama Murthy, "Efficient Transient Analysis of Finite State Space Continuous Time Markov Chains (CTMCs) : Signal Processing Approach," IEEE CACS-2013, Nantou, Taiwan